\documentclass[11pt]{article}
\usepackage{url}
\usepackage{pdfsync}
\usepackage{amssymb}
\usepackage{amsmath}
\usepackage{amsthm}
\usepackage{multirow}
\usepackage{graphicx, float}
\usepackage{url}
\usepackage{hyperref}
\usepackage[left=1in,top=1in,right=1in,bottom=1in,nohead]{geometry}
\newtheorem{theorem}{Theorem}[section]
\begin{document}
\title{Optimal airline de-ice scheduling}
\author{Jakob Kotas\thanks{Corresponding author: Department of Mathematics, University of Portland, Portland OR 97203 USA and Department of Mathematics, University of Hawai`i at Manoa, Honolulu HI 96822 USA; kotas@up.edu}, Andrew Bracken\thanks{Horizon Air, Portland OR 97218 USA; andrewbracken@gmail.com}}
\maketitle
\begin{abstract}
We present a decision support framework for optimal flight rescheduling on an airline's day of operations when de-icing  becomes necessary due to snow and ice events. Winter weather, especially in areas where such weather is not commonplace, often causes cascading delays and cancellations throughout the system due to the unforeseen need to add de-ice time to each aircraft's turnaround time. Our model optimally reschedules remaining flights of the day to minimize system delays and cancellations. The model is formulated as a mixed integer linear program (MILP). Structural properties of the model allow it to be decomposed into a finite set of linear programs (LP) and a computationally tractable algorithm for its solution is described. Finally, numerical simulations are presented for a case study of Horizon Air, a regional airline based in the Pacific Northwest of the United States.
\end{abstract}


\section{Introduction}
Operating an airline is a logistical challenge. On a normal day, careful planning months in advance brings airplanes, airline personnel, and passengers to the right airport at the right time. Small operational disruptions are accounted for and create minimal impacts. However, large events quickly overwhelm the system and bring the entire operation to a stop. Snow and ice are one such event which can become paralyzing even for the most well-prepared airlines.

Snow and ice must be removed prior to takeoff from the wings and movable control surfaces.  Without de-icing, control surfaces can become jammed and the additional drag on the wings can lead to a crash.  While there are a variety of methods used for de-icing, a clear majority of airlines use de-icing fluids which are classified based on two main characteristics: the fluid's Lowest Operational Use Temperature (LOUT) and the fluid's Holdover Time (HOT). A flight will be de-iced using fluids chosen based on the temperature (LOUT) and the estimated waiting time before takeoff (HOT). If the LOUT and/or the HOT are exceeded prior to takeoff, the flight must be de-iced again.

At airports where de-icing is a regular event, flights will depart the gate and move to a dedicated de-icing pad. This helps to streamline airport operations by allowing gates to be used by incoming flights. Additionally, resources are better allocated by having all the de-icing trucks, de-icing personnel, and airplanes collocated. At airports where de-icing is less frequent, it is done at the gate.  In this scenario, resources are spread out. De-icing trucks must drive around to airplanes across the airport. Ultimately, resources are limited to such an extent that an airport cannot support a full schedule for flights.

Deciding which flights to de-ice and in which order to de-ice them becomes a necessary but exceedingly difficult task. A flight waiting at the gate to be de-iced impacts the ability for that gate to be used by an arriving flight. Passengers waiting on landed arriving flights will miss connections.  Flight attendants and pilots that need to change flights to continue their work day will be displaced.  Delaying and cancelling flights can help to ease the burden on the system, but making the optimal decision is not obvious. During snow and ice events, airlines frequently find themselves with airplanes, airline personnel, and passengers in the wrong places at the wrong times.

In this paper we present a mathematical model for optimizing the re-scheduling of remaining flights for the day when winter weather begins. The model minimizes both the number of cancellations and the total flight delays in the system. We build up the model as a mixed-integer linear program (MILP). We then show that structural properties of the model allow it to be decomposed into a finite set of linear programs (LP). The constraints account for increased turn-around time at the gate for de-icing. This model is most appropriate for airports in which winter weather is uncommon and thus preliminary schedules do not take de-icing into account, and where de-icing is carried out at the gate for each individual aircraft.

\section{Literature review}\label{sec:litreview}
Airline operations has been an area of interest in the field of operations research (OR) for decades. The determination of timing of each flight (flight scheduling) is only one piece of the airline scheduling problem; other subproblems include crew scheduling, fleet assignment, and aircraft routing.\cite{barnhart} The full airline scheduling problem is regarded as computationally intractable. The traditional approach is to decompose the scheduling problem into its separate subproblems which are solved separately (yet suboptimally). However, some progress has been made on integrated approaches that handle multiple stages of the entire airline scheduling problem at once.\cite{papadakos,gao,weide}

Our work fits under a broad body of work devoted to recovery from system disruptions. Clarke was one of the first to give an overviews of practice in control centers under system irregularities.\cite{clarke} More recently Kohl et al. have presented an introduction to disruption management practice in industry.\cite{kohl} Within the operations research literature, Clausen et al. have presented an overview of commonly used network models for disruption recovery with references to many of the existing models for aircraft recovery.\cite{clausen} Existing formulations of the aircraft recovery problem have frequently used integer programming (IP) and mixed-integer programming (MIP).\cite{lan,abdelghany,ahmadbeygi} Some authors have also considered formulation as a minimum cost network flow problem.\cite{jarrah,mathaisel,yan,cao} Liu et al. used a multi-objective genetic algorithm for schedule disruption recovery for short-haul flights.\cite{liu}

While there has been a large amount of literature on disruptions in general, fewer authors have looked at winter weather disruptions specifically. Snow and ice disruptions are unique in that delays mostly accumulate due to the need for de-icing aircraft, which, as part of aircraft operations, can be modeled mathematically.

Norin et al. developed a heuristic framework for routing of de-ice trucks within a single airport to minimize a combination of de-ice truck distance travelled and aircraft delays.\cite{norin} Janic used deterministic queueing models to predict how snowfall creates delays due to reduced service rate of runways and gate availability; the costs of associated flight re-routes and cancellations are estimated.\cite{janic} Mao et al. described a heuristic for multi-agent-based de-ice scheduling, where the decision of which de-icing time slot to choose is made by multiple parties instead of being centrally planned.\cite{mao} While each of these papers does consider optimizing some aspect of the de-ice procedure in airline operations, we are unaware of any literature that considers optimal flight re-scheduling due to the unplanned necessity of de-icing in the way that we envision here.

\section{Model with no cancellations}\label{sec:model_wo_cancel}
We assume that decision-making based on snow events happens with very little lead time, so that a complete schedule for the day has already been developed. In particular, we assume that the departure and arrival times, and origin-destination pairs for each flight are given; we also assume that a specific aircraft (also known as ``tail number") has been assigned to each flight. In practice, this is virtually always the case for major commercial airlines, as timetables are published months in advance, and aircraft assignment days in advance, whereas the decision to delay or cancel flights due to winter weather is only made minutes to hours in advance. In this section we model the optimal re-scheduling for all flights after the ``snow-on button" has been pressed at any particular airport, without allowing cancellations. The assumption not to allow cancellations will be relaxed in section \ref{sec:model_w_cancel}. The problem is modeled with a linear program (LP). The objective function to be minimized is the weighted sum of delays to each flight in the system. We assume that the airline operates via the hub-and-spoke, rather than point-to-point, flight system, because the interconnectedness of these systems make cascading delays and cancellations more prevalent. In practice, most major airlines in the US and EU, with the notable exception of some low-cost carriers, operate via hub-and-spoke.

Let $N$ be the set of flights for the day in the entire system and let $n$ be the size of $N$. Let flights be numbered $1,2,...,n$ in such a way that flights are ordered sequentially by aircraft. For example, if aircraft \#1 operates $m_1$ flights, they are numbered $1,2,...,m_1$ in the temporal order in which they are flown, then aircraft \#2 operates $m_2$ flights which are numbered $m_1+1, m_1+2,..., m_1+m_2$, and so on. Let $s_i \in \mathbb{R}$ be the original scheduled departure time of flight $i$, for all $i=1,2,...,n$, as published in the timetable, and $x_i \in\mathbb{R}$ be the new scheduled departure time of flight $i$ after pressing the snow-on button. The delay encountered by flight $i$ is then $x_i-s_i$. All times are given in minutes after the start of the day's operations, which is typically in the early morning, in a fixed reference time zone.

We penalize delays through introduction of a cost function. We assume that costs due to delays are additive across different flights. We assume that the cost assigned to the delay encountered by flight $i$ is proportional to the delay time in minutes with a constant of proportionality $w_i\in\mathbb{R_+}$. In practice, this weighting could be equal across aircraft, or proportional to the number of passengers on the aircraft, or some other non-negative weight. Let $s\in\mathbb{R}^n$, $x\in\mathbb{R}^n$, and $w\in\mathbb{R}^n$ be column vectors whose components are $s_i$, $x_i$, and $w_i$, for $i=1,2,...,n$. $x$ is then the decision variable whereas $s$ and $w$ are known constants. The objective function is then $\min_x w(x-s)'$ where $'$ denotes transpose.

We now discuss constraints of the LP. The first constraint enforces that no flight may be scheduled earlier than originally planned in the timetable. This is reasonable as passengers, crew, and ground staff in general will not be ready to board any flight early without advance warning. Thus, $x\geq s$.

The second constraint enforces that no flight may depart before the beginning of the day's operations, in the time zone of the departing flight. Here we must convert local time to the fixed reference time. Let $z_{o,i}\in\mathbb{R}$ be the offset of the time zone of the origin airport with respect to the reference time for flight $i$, and let $z_{o}\in\mathbb{R}^n$ be the vector of $z_{o,i}$ for all $i=1,2,...n$. In other words if Pacific Standard Time is set as the reference time, then $z_{o,i}$ for a flight departing from Seattle at the beginning of its day is 0 whereas $z_{o,i}$ for a flight departing from New York at the beginning of its day, which is 3 time zones ahead (east) is 180. Then $x\geq -z_o$. In other words, if we define the beginning of the day as 5AM, then flights departing New York can leave as early as 2AM PST ($z_{o,i}=-180$) whereas flights departing Seattle can leave as early as 5AM PST ($z_{o,i}=0$).

The third constraint enforces that no flight may arrive after the end of the day's operations, which varies depending on whether aircraft are scheduled to undergo overnight maintenance procedures. Let $z_{d,i}\in\mathbb{R}$ be the offset of the time zone of the destination airport with respect to the reference time for flight $i$. Let $r_i\in\mathbb{R}_{++}$ be the scheduled duration of flight $i$ in minutes. Let $t_i\in\mathbb{R}_{++}$ be the minimum turnaround time, in minutes, that an aircraft must be on the ground after arriving before departing for its next flight. Let $d_i\in\mathbb{R}_+$ be the time to de-ice the aircraft on flight $i$ before departing. Let $e_i\in\mathbb{R}$ be the time of the end of the day's operations in the local time zone of the destination airport for flight $i$. If no maintenance is scheduled, then $e_i=24*60=1440$. Let $z_d\in\mathbb{R}^n$, $r\in\mathbb{R}^n$, $t\in\mathbb{R}^n$, $d\in\mathbb{R}^n$, and $e\in\mathbb{R}^n$ be the vectors of $z_{d,i}$, $r_i$, $t_i$, $d_i$, and $e_i$ for all $i=1,2,...,n$, respectively. Then $x\leq e-z_d - r - t - d$.

The fourth and final constraint enforces the ordering of flights on each particular aircraft. For example, if a flight is scheduled to begin the day in Seattle, then fly to Portland, then New York City, then the flight from Seattle to Portland must occur before the flight from Portland to New York City. Furthermore, there is a gap between when the first and second flights may occur which is equal to the minimum turnaround time plus de-icing time. We define a flight to be a ``sunrise flight" if it is the first flight of the day of operations for a particular aircraft. Let $S$ be the set of all sunrise flights and let $\sigma$ be the size of $S$. Then, for all flights $i\in N\backslash S$, $x_i \geq x_{i-1} + r_{i-1} + t_{i-1} + d_{i-1}$. This can be written more efficiently in matrix notation. Let $r^\dagger\in\mathbb{R}^{n-\sigma}$, $t^\dagger\in\mathbb{R}^{n-\sigma}$, and $d^\dagger\in\mathbb{R}^{n-\sigma}$ be the vectors $r$, $t$, and $d$ with all elements $i\in S$ removed. Let $M\in\mathbb{R}^{(n-\sigma)\times n}$ be a matrix that is assembled in the following way. Take the banded matrix $B\in\mathbb{R}^{n\times n}$ whose main diagonal is -1 and first upper diagonal is $1$. Then, remove all rows in $B$ corresponding to $i\in S$; this is matrix $M$. The fourth constraint can then be written as $Mx \geq r^\dagger+t^\dagger+d^\dagger$.

\begin{table}[h]
\begin{center}
{\begin{tabular}{|c|c|}
\hline
Variable&Description\\
\hline
$N$&Set of flights for the day, ordered sequentially by aircraft\\
$n$&Number of flights for the day\\
$m_i$&Number of flights flown by $i$th aircraft\\
$s$&Vector of original scheduled depart times\\
$x$&Vector of new scheduled depart times\\
$w$&Vector of weights for delay minutes on each aircraft\\
$z_o$&Vector of time zones of origin airports\\
$z_d$&Vector of time zones of destination airports\\
$e$&Vector of time of end of day's operations, minus maintenance activities\\
$r$&Vector of flight durations\\
$t$&Vector of turnaround times\\
$d$&Vector of de-ice delays\\
$S$&Set of sunrise flights\\
$\sigma$&Number of sunrise flights\\
$r^\dagger$&Vector of flight durations for non-sunrise flights\\
$t^\dagger$&Vector of turnaround times for non-sunrise flights\\
$d^\dagger$&Vector of de-ice delays for non-sunrise flights\\
$B$&Banded matrix with -1 on main diagonal and +1 on first upper diagonal\\
$M$&Matrix with rows of $B$ corresponding to sunrise flights removed\\
\hline
\end{tabular}}
\caption{Table of variables and parameters for model with no cancellations.}\label{tbl:nocancel_vars}
\end{center}
\end{table}

Thus, the LP to be solved is:
\begin{equation}\label{lp_wo_cancel_objfn}
\min_x w(x-s)'
\end{equation}
subject to
\begin{equation}\label{lp_w_cancel_constraint_1}
x\geq s
\end{equation}
\begin{equation}\label{lp_w_cancel_constraint_2}
x\geq -z_o
\end{equation}
\begin{equation}\label{lp_w_cancel_constraint_3}
x\leq e-z_d-r-t-d
\end{equation}
\begin{equation}\label{lp_wo_cancel_constraint_4}
Mx\geq r^\dagger+t^\dagger+d^\dagger
\end{equation}

\section{Cancellations}\label{sec:model_w_cancel}
\subsection{Assumptions}
Allowing flights to be canceled significantly increases the difficulty of the problem. This is due to the fact that in general, the full fleet assignment problem, where an aircraft is assigned to each particular flight, must be re-solved. For example, consider an aircraft whose scheduled flights for the day are A$\xrightarrow[]{1}$ B $\xrightarrow[]{2}$ A $\xrightarrow[]{3}$ C $\xrightarrow[]{4}$ B, where numbers above arrows denote flight numbers, and A, B, C are distinct airports. In isolation, canceling any single flight causes an issue where the aircraft is not in the correct place to carry out later flights. Even the cancellation of flight 4, the last flight of the day, does not bring the aircraft in the correct position to begin the next day's flights. There may be another aircraft in the system that can change its own locations so as to have the correct aircraft in the correct cities in the correct order. However, the full problem is very difficult, and overshadowed by the need for a quick solution in the case of a sudden winter weather event. For this reason, we make certain simplifying assumptions to make the cancellation problem tractable.

One very restrictive assumption would be to only allow cancellations to happen in pairs where the origin of the first flight and destination of the subsequent flight are equivalent: A$\rightarrow$B$\rightarrow$A. However, this carries its own issues. In many such cases, A is a hub for the airline whereas B is not; in this case, passengers going to or from location B have to wait for a different aircraft in the system, which may be many hours later, or not occur again in the same day of operations. Furthermore, pairs A$\rightarrow$B$\rightarrow$A are not always very common, especially for larger airlines with multiple hubs, where A$\rightarrow$B$\rightarrow$C and A and C being distinct hubs occurs more frequently.

Instead, we make the less restrictive assumption that flights may only be canceled if they are from one hub to another. Then, one flight either before or after the canceled flight on the same aircraft must be re-routed to adjust its origin or destination accordingly. To illustrate this, let H$_1$ and H$_2$ be two distinct hubs, and B be an airport that is not a hub. Consider an aircraft scheduled to perform B$\xrightarrow[]{1}$H$_1$$\xrightarrow[]{2}$H$_2$. Flight 2 can be canceled if flight 1 is reassigned to B$\xrightarrow[]{1}$H$_2$. Because H$_1$ and H$_2$ are both hubs, we assume there are many more flights per day between them compared to flights departing from B to any hub. Thus, passengers departing from B are inconvenienced in having to take a second connecting flight, but this is preferable to being stuck in B with no flight out whatsoever. A similar issue on an aircraft scheduled to perform H$_1$$\xrightarrow[]{1}$H$_2$$\xrightarrow[]{2}$B can be resolved by canceling flight 1 and reassigning the origin of flight 2 to H$_1\xrightarrow[]{2}$B. We also must take care not to cancel too many flights between H$_1$ and H$_2$, or else there is no added benefit for keeping passengers waiting at a hub rather than a non-hub. However, in practice the number of cancellations is small compared to the number of flights between hubs.

The incentive for this less restrictive assumption comes from the practice of Horizon Air, which we will discuss in the numerical simulation of section \ref{sec:horiz}. In looking at which flights were canceled on a day of winter weather, only flights between hubs were canceled, with adjacent flights being re-routed, as described above. For this reason we believe our assumption to be realistic and in line with current industry practice.

Another benefit of only considering flights between hubs is that airlines typically keep any and all spare planes at hubs. Thus, an aircraft that was not even scheduled for the day of operations, or was only scheduled to fly for a portion of the day of operations, can be called into duty to accommodate any necessary movement between hubs of crew, excess passengers waiting at one hub, etc.

Finally, we mention that this assumption is the most practical when H$_1$ and H$_2$ are geographically nearby hubs for two reasons: first, there are generally more flights per day between nearby hubs than far-away hubs, and second, it requires less time for passengers and crew to reposition between hubs when a cancellation does occur. For Horizon Air, which we discuss at length in section \ref{sec:horiz}, the two major hubs are Seattle and Portland, a short flight apart.

\subsection{Model with cancellations}
Let $C$ be the set of candidate flights: flights whose origin and destination airports are hubs which are near to each other, and let $c$ be the size of $C$.

First assume that we have an arbitrary set of flights that we know we want to cancel. Let the set of the flights to be canceled be $\Gamma \subseteq C$.

We model the re-scheduling problem with known cancellations $\Gamma$ using a mixed-integer linear program (MILP) which is an extension of the model of section \ref{sec:model_wo_cancel}. We introduce a cancellation penalty $p_i$ associated with flight $i$ for all $i\in C$ and let $p\in\mathbb{R}^c$ be the column vector of $p_i$ for $i\in C$. Let $y\in\mathbb{R}^c$ be a boolean vector where $y_i=1$ if $i\in\Gamma$ and $y_i=0$ if $i\notin\Gamma$, for all $i\in C$, and let $y\in\mathbb{R}^c$ be the column vector of $y_i$ for $i\in C$.

We adjust the objective function of the LP to minimize a weighted sum of delays and cancellation penalties, as $\displaystyle \min_{x,y} w(x-s)' + py'$. The relative values of $w$ and $p$ determine the decision-maker's comfort with canceling flights. For example if we take $w$ to be a column vector of all 1s of length $n$ so that every flight is weighed equally in terms of delay, then $p_i$ represents the penalty associated with canceling flight $i$ in units of effective minutes of delay per flight.

The constraints are adjusted so that the duration, turnaround, and de-ice delays associated with the canceled flight are set to zero. Thus, the canceled flight is not eliminated from the algorithm but instead ignored, as the subsequent flight can be scheduled as early as the departure time of the canceled flight. Let $r_{\Gamma,i}=r_i$ if $i\notin\Gamma$ and $r_{\Gamma,i}=0$ if $i\in\Gamma$. This sets the duration of flight $i$ effectively to 0 if flight $i$ is canceled. Define $t_{\Gamma,i}$, $d_{\Gamma,i}$, $r^\dagger_{\Gamma,i}$, $t^\dagger_{\Gamma,i}$, and $d^\dagger_{\Gamma,i}$ in the same way. Then define $r_\Gamma\in\mathbb{R}^n$, $t_\Gamma\in\mathbb{R}^n$, $d_\Gamma\in\mathbb{R}^n$ to be the vectors of $r_{\Gamma,i}$, $t_{\Gamma,i}$, and $d_{\Gamma,i}$ for $i\in N$ respectively, as well as $r_\Gamma^\dagger\in\mathbb{R}^{n-\sigma}$, $t_\Gamma^\dagger\in\mathbb{R}^{n-\sigma}$, $d_\Gamma^\dagger\in\mathbb{R}^{n-\sigma}$ to be the vectors of $r^\dagger_{\Gamma,i}$, $t^\dagger_{\Gamma,i}$, and $d^\dagger_{\Gamma,i}$ for $i\in N\backslash S$ respectively.

\begin{table}[h]
\begin{center}
{\begin{tabular}{|c|c|}
\hline
Variable&Description\\
\hline
$C$&Set of cancellable flights\\
$c$&Number of cancellable flights\\
$\Gamma$&Set of cancelled flights\\
$p$&Vector of cancellation penalties\\
$y$&Boolean vector of flight cancellations\\
$r_\Gamma$&Vector of flight durations with flights in $\Gamma$ set to 0\\
$t_\Gamma$&Vector of turnaround times with flights in $\Gamma$ set to 0\\
$d_\Gamma$&Vector of de-ice delays with flighrts in $\Gamma$ set to 0\\
$r_\Gamma^\dagger$&Vector of non-sunrise flight durations with flights in $\Gamma$ set to 0\\
$t_\Gamma^\dagger$&Vector of non-sunrise turnaround times with flights in $\Gamma$ set to 0\\
$d_\Gamma^\dagger$&Vector of non-sunrise de-ice delays with flights in $\Gamma$ set to 0\\
\hline
\end{tabular}}
\caption{Table of additional variables and parameters for model with cancellations.}\label{tbl:cancel_vars}
\end{center}
\end{table}

Since we take $\Gamma$ as given, $y$ is not a decision variable but a constant. Also note that $y$ does not appear in the constraints directly. Thus the problem for known $\Gamma$ is actually another LP:

\begin{equation}\label{lp_w_cancel_objfn}
\min_{x} w(x-s)' + py'
\end{equation}
subject to
\begin{equation}\label{lp_w_cancel_constraint_1}
x\geq s
\end{equation}
\begin{equation}\label{lp_w_cancel_constraint_2}
x\geq -z_o
\end{equation}
\begin{equation}\label{lp_w_cancel_constraint_3}
x\leq e-z_d-r_\Gamma-t_\Gamma-d_\Gamma
\end{equation}
\begin{equation}\label{lp_w_cancel_constraint_4}
Mx\geq r^\dagger_\Gamma+t^\dagger_\Gamma+d^\dagger_\Gamma
\end{equation}

However, since $\Gamma$ is not known, we must iterate this LP over every possible $\Gamma\subseteq C$. Thus, the overall problem to be solved is:

\begin{equation}\label{fullproblem_w_cancel_objfn}
\min_{\Gamma\subseteq C} \Big( \min_{x} w(x-s)' + py' \Big)
\end{equation}

where each inner problem is subject to the constraints
\begin{equation}\label{fullproblem_w_cancel_constraint_1}
x\geq s
\end{equation}
\begin{equation}
x\geq -z_o
\end{equation}
\begin{equation}
x\leq e-z_d-r_\Gamma-t_\Gamma-d_\Gamma
\end{equation}
\begin{equation}\label{fullproblem_w_cancel_constraint_4}
Mx\geq r^\dagger_\Gamma+t^\dagger_\Gamma+d^\dagger_\Gamma
\end{equation}

and $y$, $r_\Gamma$, $t_\Gamma$, $d_\Gamma$, $r_\Gamma^\dagger$, $t_\Gamma^\dagger$, and $d_\Gamma^\dagger$ are functions of $\Gamma$. We have thus replaced the MILP with an optimization over finitely many LPs by removing appropriate entries from the vectors $r$, $t$, $d$, $r^\dagger$, $t^\dagger$, and $d^\dagger$ for each inner problem.

Exhaustive search over all possible sets $\Gamma$ would thus require the solution of $2^c$ linear programs as each element of $C$ could either be canceled, or not. The exponential growth of this problem makes it computationally intractable for even modestly sized sets $C$. However, it can be proven that a much more efficient algorithm, which requires the solution of only $c+1$ linear programs, is optimal. We present the algorithm next.

\subsection{Algorithm}
Pseudo-code for the algorithm is presented below.

\begin{itemize}
\item Solve the linear program of equations (\ref{lp_w_cancel_objfn}) to (\ref{lp_w_cancel_constraint_4}) with $\Gamma = \emptyset$. This is equivalent to the linear program of equations (\ref{lp_wo_cancel_objfn}) to (\ref{lp_wo_cancel_constraint_4}). Let an optimal value of $x$ be $x^*$ and let the objective function be $f$ so that the optimal objective function value is $f(x^*)$.
\item For every $i\in C$, solve the LP of equations (\ref{lp_w_cancel_objfn}) to (\ref{lp_w_cancel_constraint_4}) with $\Gamma=\{i\}$. Let an optimal value of $x$ in this LP be $x_i^*$.
\item Let $\Gamma^*$ be the set of all $i$ for which $f(x_i^*)<f(x^*)$.
\end{itemize}

\begin{theorem}\label{thm:gammaopt}
$\Gamma^*$ is the optimum over all sets $\Gamma\subseteq C$ for the outer problem of equation (\ref{fullproblem_w_cancel_objfn}).
\end{theorem}

\subsection{Proof of theorem \ref{thm:gammaopt}}

Let $\Gamma_0\subset C$ be an arbitrary strict subset of $C$, and let $\Gamma_\gamma = \Gamma_0 \cup \{\gamma\}$ where $\gamma\in C$ and $\gamma\notin\Gamma_0$. Thus by definition $\Gamma_0 \subset \Gamma_\gamma \subseteq C$. Let $R_0$ be the feasible region defined by constraints (\ref{lp_w_cancel_constraint_1}) - (\ref{lp_w_cancel_constraint_4}) for $\Gamma=\Gamma_0$ and let $R_\gamma$ be the feasible region defined by constraints (\ref{lp_w_cancel_constraint_1}) - (\ref{lp_w_cancel_constraint_4}) for $\Gamma=\Gamma_\gamma$. We first show that $R_0 \subset R_\gamma$. Constraints (\ref{lp_w_cancel_constraint_1}) and (\ref{lp_w_cancel_constraint_2}) do not depend on $\Gamma$, so we focus on constraints (\ref{lp_w_cancel_constraint_3}) and (\ref{lp_w_cancel_constraint_4}).

First, constraint (\ref{lp_w_cancel_constraint_3}). Let $v_0 = e - z_d - r_{\Gamma_0} - t_{\Gamma_0} - d_{\Gamma_0}$ and $v_\gamma = e - z_d - r_{\Gamma_\gamma} - t_{\Gamma_\gamma} - d_{\Gamma_\gamma}$. Let $S_0=\{x|x\leq v_0\}$ and $S_\gamma=\{x|x\leq v_\gamma\}$. Recall that $r_\Gamma=r$ in all components except for those corresponding to canceled flights, in which case those components of $r_\Gamma$ are zero. Also, all components of $r>0$. Thus, $r_{\Gamma_\gamma}=r_{\Gamma_0}$ in all components except for the $\gamma^{\text{th}}$ component, in which $r_{\Gamma_0} > 0 = r_{\Gamma_\gamma}$. The vector $t > 0$ follows the same argument. The vector $d \geq 0$, which results in $d_{\Gamma_\gamma}=d_{\Gamma_0}$ in all components except for the $\gamma^{\text{th}}$ component, in which $d_{\Gamma_0} \geq 0 = d_{\Gamma_\gamma}$. Then $v_{\Gamma_0}=v_{\Gamma_\gamma}$ in all components except the $\gamma^{\text{th}}$ component, in which case $v_{\Gamma_0} < v_{\Gamma_\gamma}$. All $x\in S_0$ must also be in $S_\gamma$ since $x\leq v_{\Gamma_0} < v_{\Gamma_\gamma}$. However there exist $x \in S_\gamma$ that are not in $x\in S_0$, for example, $x=(v_{\Gamma_0} + v_{\Gamma_\gamma})/2$. Thus $S_0 \subset S_\gamma$.

Constraint (\ref{lp_w_cancel_constraint_4}) follows a similar argument; however, it is weaker, since $r^\dagger_\Gamma = r^\dagger_0$ in all components except for the $\gamma^{\text{th}}$ component {\it only if} $\gamma$ is not a sunrise flight. If $\gamma$ is a sunrise flight, then $r^\dagger_\Gamma=r^\dagger_0$ in all components. Let $T_0 = \{x | Mx\geq r^\dagger_{\Gamma_0}+t^\dagger_{\Gamma_0}+d^\dagger_{\Gamma_0} \}$ and $T_\gamma = \{x | Mx\geq r^\dagger_{\Gamma_\gamma}+t^\dagger_{\Gamma_\gamma}+d^\dagger_{\Gamma_\gamma} \}$. Then, following the reasoning as constraint \ref{lp_w_cancel_constraint_3}, we have $T_0 \subseteq T_\gamma$. The full derivation is omitted for brevity.

By definition, $R_0 = \{x | x\geq s\} \cap \{x | x \geq -z_o\} \cap S_0 \cap T_0$ and $R_\gamma = \{x | x\geq s\} \cap \{x | x \geq -z_o\} \cap S_\gamma \cap T_\gamma$. It follows from $S_0 \subset S_\gamma$ and $T_0 \subseteq T_\gamma$ that $R_0 \subset R_\gamma$.

Let us define LP$_\Gamma$ to be the LP defined in equations (\ref{lp_w_cancel_objfn}) - (\ref{lp_w_cancel_constraint_4}) with cancellations given by the set $\Gamma$. Let $x_{\Gamma}^*$ be an optimum value of $x$ for LP$_\Gamma$. Now let us compare LP$_{\Gamma_\gamma}$ and LP$_{\Gamma_0}$. The objective functions for both are identical up to a constant: LP$_{\Gamma_\gamma}$'s objective function is $p_\gamma$ larger than LP$_{\Gamma_0}$'s, due to the penalty of canceling flight $\gamma$ over the cancellations described by $\Gamma_0$. Since both objective functions have the same gradient, and $R_0 \subset R_\gamma$, there are two possibilities: either $x^*_{\Gamma_\gamma} = x^*_{\Gamma_0}$ or $x^*_{\Gamma_\gamma} \neq x^*_{\Gamma_0}$.

{\bf (a)} If $x^*_{\Gamma_\gamma} = x^*_{\Gamma_0}$, this implies that $x^*_{\Gamma_\gamma} \in R_0$ since $x^*_{\Gamma_0} \in R_0$. Thus $x^*_{\Gamma_\gamma} \notin R_\gamma \backslash R_0$, or in other words, the feasible region $R_0$ was extended to $R_\gamma$ in (a) direction(s) that did not improve the optimal objective function value. Since the objective function for LP$_{\Gamma_\gamma}$ at $x^*_{\Gamma_0}=x^*_{\Gamma_\gamma}$ is greater than the objective function of the LP$_{\Gamma_0}$ at the same point by an amount $p_\gamma>0$, when considering the full problem \ref{fullproblem_w_cancel_objfn} - \ref{fullproblem_w_cancel_constraint_4}, $\Gamma_0$ would be more optimal than $\Gamma_\gamma$ for the outer problem. Therefore, any time we try to cancel a new flight $\gamma$ and the $x^*_{\Gamma_\gamma} = x^*_{\Gamma_0}$, we may immediately conclude that canceling flight $\gamma$ is never optimal.

{\bf (b)} If $x^*_{\Gamma_\gamma} \neq x^*_{\Gamma_0}$, this implies that $x^*_{\Gamma_\gamma} \in R_\gamma \backslash R_0$. $x^*_{\Gamma_\gamma} \in R_\gamma$ because it is the solution to LP$_{\Gamma_\gamma}$, but it is not in $R_0$ because $wx^*_{\Gamma_\gamma} < wx^*_{\Gamma_0} \leq wx_{\Gamma_0}$ for all $x_{\Gamma_0} \in R_0$. We must then compare the objective function values of LP$_{\Gamma_0}$ at $x^*_{\Gamma_0}$ and LP$_{\Gamma_\gamma}$ at $x^*_{\Gamma_\gamma}$. If $w{x^*}_{\Gamma_0}' + p_{\gamma} < w{x^*}_{\Gamma_\gamma}'$, then canceling flight $\gamma$ has improved the objective function value as compared to not having canceled flight $\gamma$.

Let $\Gamma^*$ be the optimal $\Gamma$ to the outer problem of equations (\ref{fullproblem_w_cancel_objfn}) - (\ref{fullproblem_w_cancel_constraint_4}). It follows from {\bf (b)} that $x^*_{\Gamma^*} \in R_{\Gamma^*} \backslash R_{\Gamma^* \backslash \{i\}}$ for all $i$. At the same time, it follows from {\bf (a)} that $x^*_{\Gamma^*}\in R_{\Gamma^*\backslash\{ i\}}$ and thus $x^*_{\Gamma^*} = x^*_{\Gamma^*\backslash\{ i\}}$. Thus for all $j\notin \Gamma^*$, any set $\Gamma \ni j$ is suboptimal compared to that same $\Gamma \cup \{ j \}$.

\begin{figure}[h]
\begin{center}
{\begin{tabular}{llll}
{\bf a.}&&&\\
&&\includegraphics[height=3.25cm]{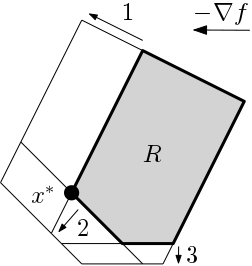}\\
{\bf b.}&{\bf i.}&{\bf ii.}&{\bf iii.}\\
&\includegraphics[height=3cm]{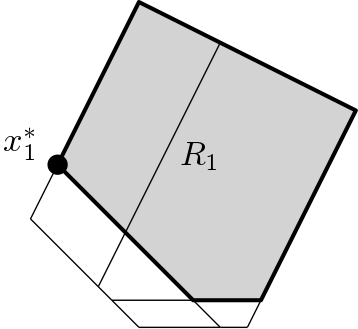}&\includegraphics[height=3cm]{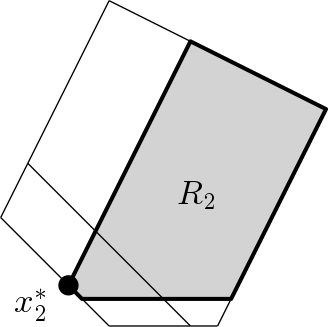}&\includegraphics[height=3cm]{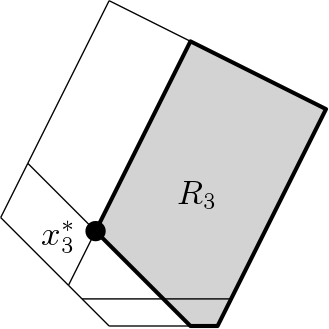}\\
{\bf c.}&{\bf i.}&{\bf ii.}&{\bf iii.}\\
&\includegraphics[height=3cm]{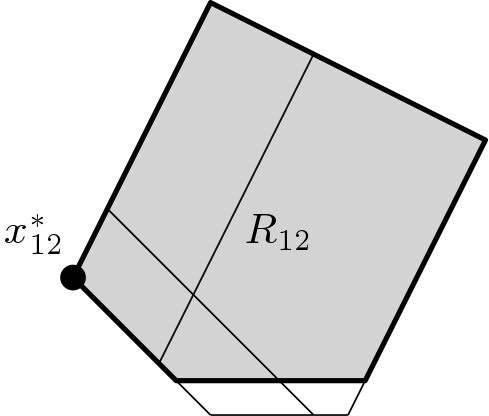}&\includegraphics[height=3cm]{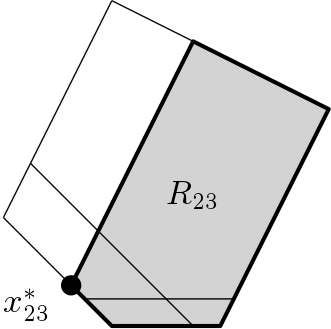}&\includegraphics[height=3cm]{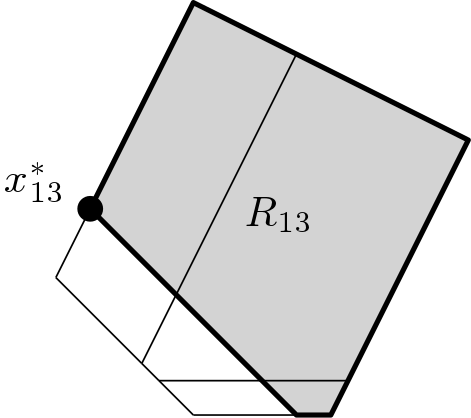}
\end{tabular}}
\caption{A visual aid to Theorem \ref{thm:gammaopt}. Subfigure (a) shows a shaded polygonal feasible region $R$ in 2D with $f(x)$ being the objective function value, left being the direction if decreasing $f(x)$ and $x^*$ being the optimal point in region $R$. Arrows with numbers denote directions in which feasible regions are extended when a certain flight is canceled. In subfigures (b) and (c), $R_{ab}$ denotes feasible regions extended when flight(s) $a$ (and $b$) are canceled, and $x^*_{ab}$ is the corresponding optimal point. In this example, $x^*_{12}$ is the farthest left of any optimum point in any subfigure. We see that canceling either flight $1$ or $2$ always brings the optimum farther to the left (for example, $x^*_2$ is to the left of $x^*$ and $x^*_{13}$ is to the left of $x^*_3$) while canceling flight $3$ does not result in movement of the optimal point (for example, $x^*_{13} = x^*_1$).}\label{fig:feas_region}
\end{center}
\end{figure}

A visual aid to Theorem \ref{thm:gammaopt} is given in figure \ref{fig:feas_region}.

\section{Model with cancellations applied to Horizon Air schedule}\label{sec:horiz}
Horizon Air is a regional airline and subsidiary of Alaska Air Group. Since 2011, it has followed a capacity purchase agreement business model where all Horizon-operated flights are marketed and sold by Alaska Airlines. Horizon focuses mainly on flights to, from, and within the Pacific Northwest region of the United States, with Seattle as its main hub and Portland a secondary hub.

Seattle and Portland experience oceanic climates with cool wet winters and mild dry summers. Rain is frequent in the winter months of October to March. Light snow is not uncommon, but heavy ice and snowstorms are rare, occurring less than once per year, on average.

On December 24-25, 2017, a snowstorm occurred in the Pacific Northwest, with snow in Seattle and light snow and freezing rain in Portland. This led to the need for de-icing flights departing Seattle and Portland, with many flights subsequently delayed, and some canceled. Data for scheduled and actual departure and arrival times for all flights flown in the US is publicly available from the Bureau of Transportation Statistics.\cite{bts} Data was obtained for the scheduled and actual departure and arrival of all flights flown by Horizon Air on Monday, December 25, 2017. Out of 276 regularly scheduled flights, 119 (43.1\%) were delayed by over 15 minutes, with 54 (19.6\%) experiencing significant delays of over one hour. 16 (5.8\%) were canceled, with 7 extra re-routing flights being added to move aircraft to the correct location. One such example of a cancellation and added flight is given in table \ref{tbl:cancel_example}. The total system-wide delays (sum of actual departures minus scheduled departures, with early departures counting as 0, and not counting canceled flights) was 7787 minutes, or an average of 29.95 minutes per (non-canceled) flight. Compared to a more typical Monday without winter weather, one week prior on December 18, 2017, just 1 out of 334 flights was canceled (0.3\%) while the total system-wide delays was 7880 minutes, or 23.66 minutes per (non-canceled) flight.

\begin{table}[h]
\begin{center}
{\begin{tabular}{|c|c|c|c|c|}
\hline
Flight & Origin & Scheduled Departure & Destination & Scheduled Arrival\\
\hline
\multicolumn{5}{|c|}{Canceled:}\\
\hline
2473 & SEA & 9:45 & PDX & 10:44\\
2209 & PDX & 11:32 & MFR & 12:29\\
\hline
\multicolumn{5}{|c|}{Added:}\\
\hline
9372 & SEA & 11:15 & MFR & 12:18\\
\hline
\end{tabular}}
\caption{An example of two flights being canceled, with one added, so that the aircraft is in place at MFR for later afternoon flights. In our model, we count this as one (net) cancellation (flight 2473) plus one re-routing leg (flight 2209 becomes flight 9372).}\label{tbl:cancel_example}
\end{center}
\end{table}

To solve the full problem of equations (\ref{fullproblem_w_cancel_objfn}) to (\ref{fullproblem_w_cancel_constraint_4}), we make some assumptions for parameter values and perform sensitivity analysis on them later.

$C$, the set of eligible flights for cancellation, is comprised of all flights between Horizon's two hubs of Seattle (SEA) and Portland (PDX) departing after the snow-on button has been pressed. As mentioned before, one adjacent flight will need to be re-routed so that the aircraft is in the correct location at the end of the cancellation+re-route pair. In some cases, the same aircraft performs two flights in a row between hubs (i.e., SEA$\rightarrow$PDX$\rightarrow$SEA or PDX$\rightarrow$SEA$\rightarrow$PDX). In these cases, the re-routing flight amounts to cancellation of the second flight. We assign $p_i=60$ for all $i\in C$ that have at least one adjacent flight on the same aircraft that is also $\in C$. In other words, the cancellation of one such flight is equivalent to a 60-minute delay in terms of overall objective function value. Since both flights are canceled, the overall addition to the objective function will be $2p_i=120$. These flights overall are ``easier" to cancel as the aircraft is in the correct location without needing to re-route another flight.

For all other flights $i\in C$, we assign a higher penalty, due to the need for re-routing. Intuitively we want $p_i$ for these flights to be at least twice the penalty for flights who have an adjacent SEA-PDX or PDX-SEA flight, due to the added inconvenience of re-routing a different flight. We choose $p_i=180$ for these flights. In other words, cancellation of one such flight is penalized equally to a 180-minute delay.

$w\in\mathbb{R}^n$ was taken to be a vector of ones so that one minute of delay on any particular flight is penalized equally to one minute of delay on any other flight.

Minimum turnaround time of an aircraft on the ground between arrival of the previous flight and departure of the following flight was taken to be 45 minutes across all flights. De-icing was assumed to add an extra 20 minutes to all flights departing SEA or PDX beginning when the respective snow-on button is pressed at each airport until the end of the day of operations.

We take the beginning of the operational day to be 5AM local time. We ignore overnight maintenance so that the end of the day is defined to be $e=1440$ in the local time zone of each arriving flight. That is, all flights planned for December 25 must arrive at their destination by 5AM on December 26, local time.

We assume the snow-on button is pressed for both Seattle and Portland at the beginning of the operational day on December 25. This is reasonable on this date since some snow had been falling since the day before. $C$ then  contains 29 flights out of the total 276 scheduled. For some of these flights which occur neither at the beginning or the end of the day, there is a choice of which adjacent flight (neither of which is SEA$\rightarrow$PDX or PDX$\rightarrow$SEA) to re-route. Our algorithm does not specify which to choose and regards both as equivalent. In practice, one would most likely choose the flight with the fewest passengers.

The problem was solved in MATLAB R2017b.\cite{matlab} Each inner LP of the full problem of equations (\ref{fullproblem_w_cancel_objfn}) - (\ref{fullproblem_w_cancel_constraint_4}) was solved using the CVX software package.\cite{cvx}

\subsection{Results}\label{sec:results}
The optimal solution of our model is to cancel only 2 flights, as shown in table \ref{tbl:cancel_baseline_cancelpenalty}. This is a reduction of 87.5\% compared to the 16 that were canceled on the actual day of operations. The total system-wide delays are 6470 minutes, a reduction of 16.9\% compared to 7787 minutes on the actual day of operations.

\begin{table}[h]
\begin{center}
{\begin{tabular}{|c|c|c|c|c|}
\hline
Flight & Origin & Scheduled Departure & Destination & Scheduled Arrival\\
\hline
2148 & PDX & 18:05 & SEA & 18:58\\
2211 & SEA & 19:34 & PDX & 20:23\\
\hline
\end{tabular}}
\caption{Optimal flights to cancel for the Horizon Air system on December 25, 2017.}
\end{center}\label{tbl:cancel_baseline_cancelpenalty}
\end{table}

\subsection{Sensitivity analysis}
\subsubsection{Snow-on time}
\begin{figure}[h]
\begin{center}
{\includegraphics[width=12cm]{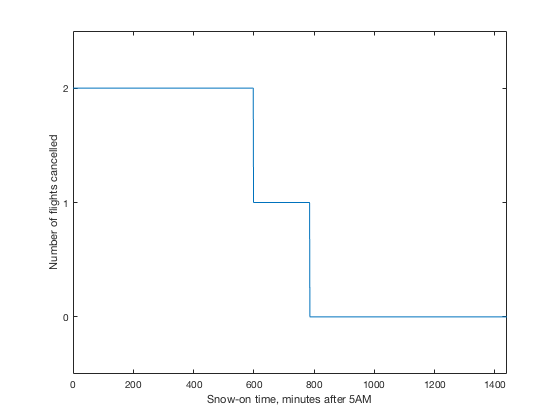}}
\caption{As snow-on time is pushed later in the day, fewer flights are canceled. For both flights, there exists a lead-in time where the flight is only canceled if the scheduled departure time is beyond a certain window after snow has begun.}\label{fig:snow_on_sa}
\end{center}
\end{figure}

We investigate the effect of changing the time at which the snow-on button is pressed at SEA and PDX. We assume that the snow-on time is the same at both hubs. We have already shown in section \ref{sec:results} that when the snow-on time is at 5AM, two flights are canceled. We now allow the snow-on time to vary over the entire day of operations, from 0 to 1440 (in minutes after 5AM.) Results are shown in figure \ref{fig:snow_on_sa}.

Flight 2148 is scheduled to depart at 18:05, or time 785. Flight 2211 is scheduled to depart at 19:34, or time 874. Numerical simulation shows that when the snow-on time is $\leq$ 598, both flights 2148 and 2211 are canceled. When the snow-on time is $\geq$ 599 and $\leq$ 785, only flight 2148 is canceled. When the snow-on time is $\geq$ 786, no flights are canceled. In both cases, we see that there exists a lead-in time; flights scheduled to depart sooner than a critical time in advance of snow are left alone, while flights beyond that critical time are canceled. (Note that this ``critical time" is not constant across flights, but is positive in both cases.) This matches common practice in industry where flights departing many hours from now would be canceled before flights departing immediately, as the increased lead-in time allows for more flexible re-scheduling of aircraft, passengers and crew.

\subsubsection{Cancellation penalty}
\begin{figure}[h!]
\begin{center}
{\includegraphics[width=12cm]{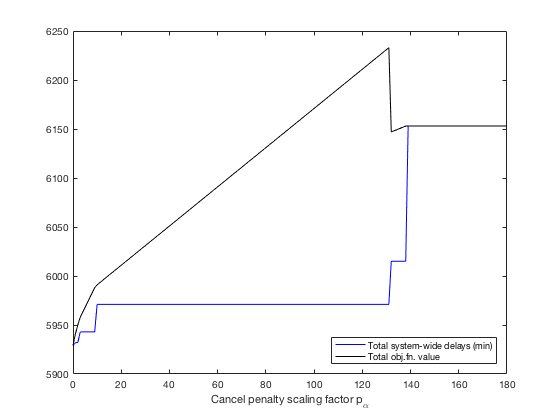}}
\caption{Total system-wide delays, in minutes, and total objective function value (delays + cancellation penalties) as a function of the cancel penalty scaling factor $p_\alpha$. As $p_\alpha$ increases, the decision-maker becomes more averse to canceling flights.}\label{fig:p_alpha_sa}
\end{center}
\end{figure}

We now consider changing the cancellation penalty vector $p$. Let cancellable flights with adjacent flights on the same aircraft that are also between the two hubs, so as to avoid the necessity of a re-routing flight, have a penalty value $p_\alpha$, and let cancellable flights with no such adjacent flight on the same aircraft, so that re-routing of an adjacent flight is necessary, have a penalty value $p_\beta$. We wish to maintain that $p_\beta > 2p_\alpha$ as discussed in section \ref{sec:horiz}. To reduce the number of degrees of freedom to one, we take $p_\beta = 3p_\alpha$. Whereas we previously took $p_\alpha=60$, we now vary $p_\alpha$ from 0 to 180. In figure \ref{fig:p_alpha_sa}, we show the total system-wide delays and total objective function value as functions of $p_\alpha$. The total objective function value is the sum of total system-wide delays plus the sum of all cancel penalties. As expected, as the cancel penalty scaling factor $p_\alpha$ increases, the decision-maker becomes more averse to canceling flights over delaying them, so the total system-wide delays increase.

\begin{figure}[h!]
\begin{center}
{\includegraphics[width=12cm]{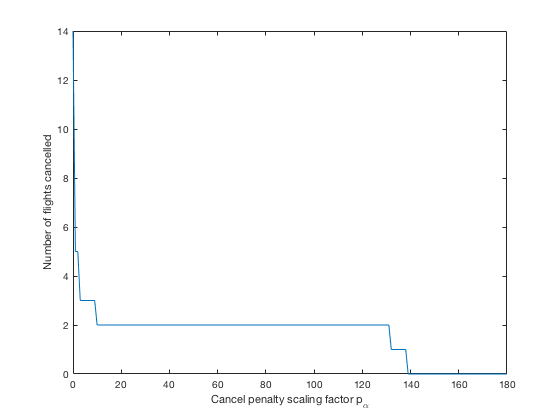}}
\caption{Number of flights canceled as a function of the cancel penalty scaling factor $p_\alpha$. As $p_\alpha$ increases, the decision-maker becomes more averse to canceling flights.}\label{fig:p_alpha_sa2}
\end{center}
\end{figure}

In figure \ref{fig:p_alpha_sa2}, we show the number of flights canceled as a function of $p_\alpha$. Again, we see the number of flights canceled decreases as the decision-maker becomes more averse to canceling flights rather than delaying them. 

In practice, a decision-maker may have an order-of-magnitude sense of a reasonable value for $p_\alpha$. However, it may be more helpful instead to consider which flights to cancel, given that we must cancel a certain number. For example, if a decision-maker is comfortable canceling as many as 4 flights, which should he or she cancel? This question can be answered by varying $p_\alpha$ from $+\infty$ to 0, keeping track of the value of $p_\alpha$ at which each flight is canceled, and ordering them. This creates a ranked list of flights to cancel, from most to least optimal. This list is given in table \ref{tbl:cancel_rank}.

We do note, however, that values of $p_\alpha$ less than perhaps 20 would indicate a decision-maker who is {\it{very}} willing to cancel flights. If $p_\alpha=20$, then canceling a flight is equivalent to only a 20- or 60-minute delay on one flight, depending on whether rerouting is needed for an adjacent flight. This is more evidence that canceling 16 flights, as was done in practice on December 25, 2017, is far from optimal, and a solution with  fewer cancellations is possible. Nevertheless we include flights in the rank who would only be canceled for very low $p_\alpha$, down to $p_\alpha=1$, in table \ref{tbl:cancel_rank} for completeness.

\begin{table}[h!]
\begin{center}
{\begin{tabular}{|c|c|c|c|c|c|c|c|}
\hline
Rank & Max $p_\alpha$ & Flight & Origin & Sched. Departure & Destination & Sched. Arrival&Canceled?\\
\hline
1&138.0&2211 & SEA & 19:34 & PDX & 20:23 & No\\
2&131.9&2148 & PDX & 18:05 & SEA & 18:58 & No\\
3&9.3&2290&PDX&9:50&SEA&10:46&Yes\\
4&2.9&2301&SEA&13:40&PDX&14:29&No\\
5&2.6&2328&PDX&23:28&SEA&00:32&Yes\\
\hline
\end{tabular}}
\caption{Ranking of best flights to cancel for the Horizon Air system on December 25, 2017. Max $p_\alpha$ indicates the value of cancel penalty scaling factor $p_\alpha$ below which canceling the designated flight is optimal. The top two flights are the only likely to be canceled for reasonably high values of $p_\alpha$. The rightmost column indicates whether or not the flight was canceled in reality on December 25, 2017. We see that only 2 of the top 5 predicted by our model were in fact canceled, as well as 14 others not shown.}\label{tbl:cancel_rank}
\end{center}
\end{table}

\section{Conclusions}
We have developed a mathematical framework for flight re-scheduling in the case of unexpected winter weather. This framework optimally readjusts flight departure times to allow extra time for aircraft de-icing so as to minimize both total system delays and cancellations. This framework is most useful for airlines with significant operations in cities that rarely (yet sometimes) experience winter weather, so that it is not taken into account in preliminary schedule formation. The model was built as a finite set of linear programs (LPs). We proved that the underlying structure of the problem allows for efficient solution over those sets. A numerical simulation was performed on data from Horizon Air, whose hubs are Seattle and Portland, on a date when winter weather impacted both of those airports. Our simulations predicted that our model would have reduced delays by 16.9\% and cancellations by 87.5\% compared to the actual day of operations. Sensitivity analysis on model parameters revealed solutions changing as expected: that pushing the snow-on time later reduces the number of cancellations as there are fewer possible flights to be canceled; and that increasing the cancellation penalty parameter decreases the number of cancellations while correspondingly increasing the total system delays.

De-icing is just a small piece of the flight scheduling puzzle, and that in itself is part of the broader problem of airline scheduling. As such, future directions for research would include incorporating our de-icing model into a more general model for disruption management due to other causes. In addition, by embedding our model within an integrated scheduling algorithm, a relaxation of our assumptions on which flights can be canceled could be developed in conjunction with a re-solving of the fleet assignment problem. Another possible area for future research would be to take flight durations, turnaround, and de-ice times to be not deterministic, but random variables with known distributions, resulting in a set of stochastic linear programs.

\section{Acknowledgments}
This work was partially funded by a Butine Grant from the University of Portland.

\bibliographystyle{plain}
\bibliography{deicing}
\end{document}